\documentclass{article}

\usepackage{xurl}

\usepackage{amsmath,amssymb,graphicx,hyperref,cleveref}
\usepackage[table]{xcolor}
\usepackage{algorithm}

\usepackage{amsthm}

\theoremstyle{plain}
\newtheorem{theorem}{Theorem}[section]

\theoremstyle{definition}

\newcommand{\tri}{\mathcal T}
\newcommand{\manifold}{\mathcal M}
\newcommand{\surface}{\mathcal S}
\newcommand{\simplex}{\Delta}
\newcommand{\subsimplex}{\delta}

\begin{document}

\title{Sampling triangulations of manifolds using Monte Carlo methods}
\author{Eduardo G. Altmann\thanks{School of Mathematics and Statistics F07, The University of Sydney, NSW 2006 Australia, \texttt{eduardo.altmann@sydney.edu.au}}
\and Jonathan Spreer\thanks{School of Mathematics and Statistics F07, The University of Sydney, NSW 2006 Australia, \texttt{jonathan.spreer@sydney.edu.au}}}

\date{}

\maketitle

\begin{abstract} \small We propose a Monte Carlo  method to efficiently find, count, and sample abstract triangulations of a given manifold $\manifold$. The method is based on a biased random walk through all possible triangulations of $\manifold$ (in the Pachner graph), constructed by combining (bi-stellar) moves with suitable chosen accept/reject probabilities (Metropolis-Hastings). Asymptotically, the method guarantees that samples of triangulations are drawn at random from a chosen probability. This enables us not only to sample (rare) triangulations of particular interest but also to estimate the (extremely small) probability of obtaining them when isomorphism types of triangulations are sampled uniformly at random. 
We implement our general method for surface triangulations and $1$-vertex triangulations of $3$-manifolds. To showcase its usefulness, we present a number of experiments: {\em (a)} we recover asymptotic growth rates for the number of isomorphism types of simplicial triangulations of the $2$-dimensional sphere; {\em (b)} we experimentally observe that the growth rate for the number of isomorphism types of $1$-vertex triangulations of the $3$-dimensional sphere appears to be singly exponential in the number of their tetrahedra; and {\em (c)} we present experimental evidence that a randomly chosen isomorphism type of $1$-vertex $n$-tetrahedra $3$-sphere triangulation, for $n$ tending to infinity, almost surely shows a fixed edge-degree distribution which decays exponentially for large degrees, but shows non-monotonic behaviour for small degrees.

\end{abstract}

%57Q15   	Triangulating manifolds
%60J10		Markov chains
%57-08		Computational methods for problems pertaining to manifolds and cell complexes
\noindent 
{\it 2020 Math. Subj. Classification.} Primary: 57Q15; Secondary: 60J10, 57-08.


\section{Introduction}
\label{sec:introduction}

The decomposition of a domain or space into simple pieces such as triangles, tetrahedra and their higher-dimensional analogues, i.e., building a {\em triangulation} of the domain or space, is a standard technique across a multitude of scientific disciplines. Triangulations are used to approximate surfaces in computer graphics, to describe domains on which partial differential equations are solved, and to efficiently represent manifolds. In most of these applications, triangulations face a fundamental challenge: they are far from unique, and the choice of a triangulation can have dramatic consequences on running times and stability of the algorithms we run on them. What is considered a good triangulation varies across applications: in discrete geometry, Delaunay triangulations are often preferred because they maximise their smallest angle between faces and achieve more stable behaviour of computations, see \cite[Section 2.2.2]{Triangulations}; in geometric topology, triangulations with the smallest number of pieces or with small structural parameters are considered good because algorithms to solve important problems run efficiently on them \cite{Burton11Simplification,Huszar19Treewidth}.  

Whenever triangulations are used to computationally solve geometric or topological problems, two related questions are relevant:
\begin{enumerate}
    \item How do we find triangulations with optimal properties (for a specific application)?
    \item What properties do we expect from a typical input triangulation? 
\end{enumerate}

Tackling these two questions in different geometric and topological settings has produced a considerable body of research, see 
\cite{BenedettiZiegler11LCSpheres, Benedetti13KnotsInSpheres, BPSNormSurfExp, Hachimori00FewEdgeKnots, Lickorish91Unshellable, Shankar_2023} for a selection, and below for more. But fundamental questions, such as the number of abstract triangulations of the $3$-dimensional sphere, remain open despite significant research efforts \cite{Kalai88ManySpheres,Pfeifle02ManySpheres,Nevo15ManySpheres}.

A variety of numerical approaches to address the questions above experimentally explores the space of triangulations using bi-stellar moves. This is what we build on here. In the geometric setting, where point configurations with coordinates are triangulated, a comprehensive overview can be found in \cite{Triangulations}. In combinatorial topology and low-dimensional geometric topology, bi-stellar moves on triangulations of manifolds have been employed to construct small triangulations of standard $3$-manifolds \cite{Bjoerner00Flips}, or examples of triangulations of spheres with interesting properties \cite{burton23flipgraph}. It has also been widely used as a pre-processing step for expensive algorithms such as $3$-sphere or unknot recognition in $3$-manifold topology \cite{Burton11Simplification,Joswig22Frontiers}. Support for bi-stellar moves is part of software in combinatorial topology \cite{simpcomp,polymake}, as well as in low-dimensional topology \cite{regina}.
A limitation of most of these computational methods is that they are unable to {\it simultaneously} tackle both questions mentioned above: the methods that successfully tackle question 1 and can find rare triangulations fail at question 2 as they do not control the distribution from which they sample so that they are unable to estimate the probabilities of triangulations (e.g., when sampled uniformly at random from the space of combinatorial isomorphism classes of triangulations).

An extremely powerful computational approach to tackle {\em both} questions above in different high-dimensional settings is to use Markov Chain Monte Carlo (MCMC) methods~\cite{Robert2004}. In particular, a variety of MCMC methods are used to sample random graphs~\cite{Cimini2019,Coolen2017}, which are based on suitably chosen moves (proposals) that preserve graph properties~\cite{Maslov2002,Milo2004} and employ different MCMC (acceptance) techniques~\cite{Byshkin2016,Fischer2015,Snijder2002}. More recently, these results have been further expanded to simplicial complexes~\cite{Young2017} as part of a general effort to generalise results from network theory to hypergraphs~\cite{Battiston2021}. Triangulations of manifolds form a particularly simple and interesting class of hypergraphs (e.g. pure simplicial complexes) and their study can contribute to the more general effort of understanding (numerical methods in) more general hypergraphs. This is in line with other successful connections between Statistical Physics and triangulations~\cite{aste2012,prd}. It is thus natural to consider MCMC methods to study triangulations, and to address the limitations of previous numerical approaches. 
 
\paragraph*{Our contribution}
In this article we propose a Monte Carlo method to find and sample abstract triangulations of manifolds (see \Cref{sec:method}). This is achieved by iteratively performing local modifications on them (using bi-stellar moves), but adapting the resulting random walk in such a way that both questions above can be answered: we can both find rare triangulations of interest and estimate probabilities of randomly selected triangulations. This is done because in the limit of large number of moves, the Markov Chain created by our method is guaranteed to visit every isomorphism type of triangulation with a pre-assigned probability. 
Our method is adaptable to arbitrary dimensions $d$, and to arbitrary triangulation types (generalised triangulations, graph-encoded manifolds, simplicial complexes, etc.). 

We implement our method in dimensions two and three (\Cref{ssec:dim2Desc,ssec:dim3Desc}), and perform several experiments to showcase its abilities.
More precisely, in \Cref{ssec:basic} we show that our method works in practice and as intended. 
In \Cref{ssec:numberoftriangulations} we confirm 
known results about the number of isomorphism classes of surface triangulations. Specifically, we reproduce \cite[Equation (8.1)]{tutte_1962} for
the number of simplicial triangulations of the $2$-sphere. We also estimate numbers of generalised triangulations of surfaces of low genus,
thereby experimentally confirming analogous asymptotic behaviour of their growth rates as in the simplicial case, see \cite{Gao91NumberOfRootedTriangulations,Goulden08KPHierarchy,Bender13ApproximationNumberTriangulations}. 
Most importantly, we give experimental evidence for the growth rate of the number of isomorphism types of $1$-vertex 
triangulations of the $3$-dimensional sphere to be only singly exponential in the number of tetrahedra. 
If the number of arbitrary $n$-tetrahedra triangulations of the $3$-sphere was super-exponential in $n$, this growth rate would be expected to show up in the subset of $1$-vertex triangulations as well.
In \Cref{ssec:vertexdegrees} we investigate the number of low-degree edges in triangulations of the $3$-sphere.
Our findings show that the average proportion of edges of degree at most five in a given triangulation converges to a fixed 
share of all of its edges. Moreover, the standard deviation of this average seems to decay in accordance with the central limit theorem. 
Thus, our experiments suggest that, for $n$ tending to infinity, a randomly chosen $n$-tetrahedron triangulation of the $3$-sphere
almost surely has a fixed edge-degree sequence, which we compute.

The above experiments are meant to highlight what can be done with our method. We believe that this method, or more efficient variations thereof, can be used to gain more insights into the space of triangulations and as a computer assisted device to find conjectures and counter-examples. Moreover, we think that extensions of the method have the potential to be applied in geometric settings, or for more general objects such as hypergraphs.

\subsection*{Acknowledgements}
This work received support from the Australian Research Council under the Discovery Project scheme (grant number DP220102588) and was finished while the authors were on sabbatical in Germany (J.S. at Technische Universit\"at Berlin and E.G.A. at the Max Planck Institute for the Physics of Complex Systems in Dresden). We thank both institutions for their support and hospitality, in particular Michael Joswig (Berlin) and Holger Kantz (Dresden).

%J.S. is supported by the Australian Research Council under the Discovery Project scheme (grant number DP220102588). This work was finished while J.S. was on sabbatical at Technische Universit\"at (TU) Berlin and E.G.A. was on sabbatical at the Max Planck Institute for the Physics of Complex Systems (MPIPKS) Dresden. J.S. wants to thank TU Berlin, and in particular Michael Joswig and his group for their hospitality.   

\section{Background: triangulations of manifolds}
\label{ssec:triangulations}

\subsection{Manifolds}
A {\em $d$-dimensional manifold}, or {\em $d$-manifold} for short, is a (topological, second-countable, Hausdorff) space $\manifold$ such that every point $p \in \manifold$ has a neighbourhood $U(p) \subset \manifold$ that is {\em homeomorphic} to, i.e, a continuous deformation of, the Euclidean space $\mathbb{R}^d$.
Manifolds in dimension $1$ are circles. $2$-manifolds are better known as surfaces, with the most prominent examples being the $2$-dimensional sphere (the surface of the Earth) and the $2$-dimensional torus (the surface of a donut). 
Famous non-examples for manifolds are the figure eight, or multiple $2$-spheres pinched together at a point.
Note that, in our definition, the disk $\{ x \in \mathbb{R}^2 \,\mid\, |x|\leq 1\}$ is not considered to be a manifold as it has boundary points.
Two manifolds $\manifold$ and $\manifold'$  are considered equivalent, or {\em homeomorphic}, if one can be transformed into the other by a continuous deformation (e.g., the surface of a coffee mug with a single handle, and the surface of a donut are considered to be the same $2$-dimensional manifold).

Often we can endow a manifold with additional structure. For instance, we can require the neighbourhoods of points $U(p)$ to be smooth or piecewise linear deformations of $\mathbb{R}^d$ (in a consistent way throughout the manifold). This leads to the notion of a {\em smooth} or {\em piecewise linear (PL) manifold}. Given two PL manifolds we say that they are {\em piecewise linearly homeomorphic}, if they can be transformed into each other by a piecewise linear motion or, equivalently, if they have a common subdivision into linear pieces.

Up to dimension three, every manifold has a unique piecewise linear structure (and a unique smooth structure), and the concepts of homeomorphy and PL homeomorphy (and smooth homeomorphy, also known as {\em diffeomorphy}) coincide. From dimension four on, there are manifolds with multiple in-equivalent PL structures, and also manifolds that do not admit any PL structures at all. See \cite{Gompf,Saveliev+2012} for more background reading on (four-)manifolds.

\subsection{Simplices}
The $d$-dimensional simplex is the convex hull of $d+1$ points in general position. In dimension $0$, $1$, $2$, and $3$, simplices are called {\em vertices}, {\em edges}, {\em triangles}, and {\em tetrahedra} respectively. We can think of a $d$-simplex $\simplex$ as a purely combinatorial object with set of points $V(\simplex) = \{ 0,1, \ldots , d \}$. Every subset of $V(\simplex)$ of cardinality $(i+1)$ spans a sub-simplex $\subsimplex \subset \simplex$ called an {\em $i$-dimensional face}, or {\em $i$-face}. A face of a simplex is called {\em proper}, if its corresponding subset is. 

\begin{figure}[htb]
  \centering{\includegraphics[width=0.7\textwidth]{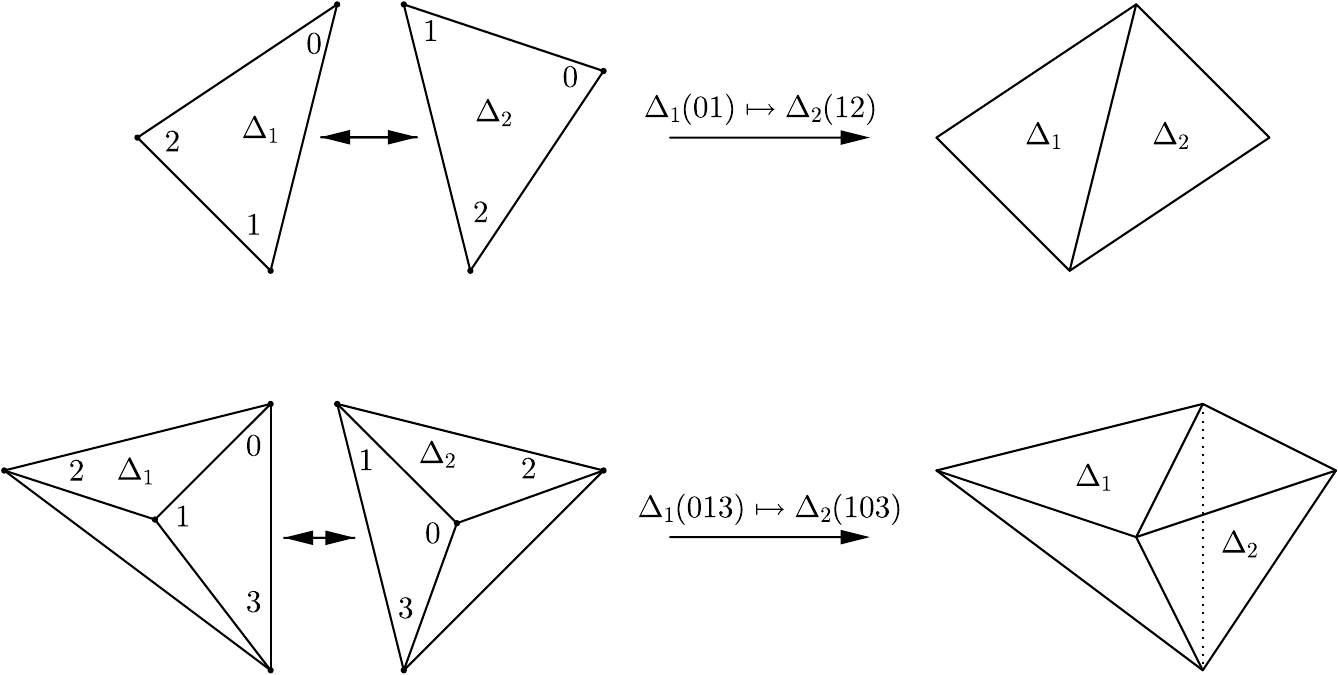} }
  \caption{Face gluings in dimensions two (top) and three (bottom). \label{fig:triangulations}}
\end{figure}

Given two $d$-simplices $\simplex_1$ and $\simplex_2$, we can glue them together along two of their $(d-1)$-dimensional faces, by specifying how $d$ of their respective $d+1$ points (or vertices) are identified with each other and then linearly interpolate. We call this a {\em face gluing} of $\simplex_1$ and $\simplex_2$, see \Cref{fig:triangulations} for examples. A face gluing can also identify a pair of $(d-1)$-dimensional faces of the same simplex.

\subsection{Triangulations of manifolds}

Informally speaking, triangulating a manifold $\manifold$ in dimension $d$ means dividing it up into $d$-simplices, referred to as {\em facets}, that are identified along face gluings. More precisely, given a set of $n$ disjoint $d$-simplices $\simplex_1, \ldots , \simplex_n$, we can glue them together along their $n \cdot (d+1)$ faces of dimensions $d-1$ in pairs. Let $\Phi_1 , \ldots , \Phi_{n \cdot (d+1)/2}$ be the corresponding face gluings. We call the quotient space
$$\tri = \{\simplex_1, \ldots , \simplex_n\} / \{\Phi_1 , \ldots , \Phi_{n \cdot (d+1)/2}\}$$ a {\em PL triangulation of a $d$-dimensional manifold}, if the following conditions are satisfied:
\begin{enumerate}
  \item If a face becomes identified with itself as a result of a face gluing, then this must happen along the identity map.
  \item The boundary of a small neighbourhood of each vertex of $\tri$, considered with its natural decomposition into $(d-1)$-simplices coming from $\tri$, must be a PL triangulation PL homeomorphic to the $(d-1)$-sphere (the set of points at unit distance in $\mathbb{R}^d$).
\end{enumerate}
The definition of a PL triangulation then becomes complete by adding that the PL triangulation of a $0$-sphere is two isolated points.\footnote{The notion of a PL triangulation of a manifold is not to be confused with the slightly different notion of a {\em triangulated manifold}. However, both terms are equivalent for dimensions at most three.} 

As a result of the face gluings, multiple lower-dimensional faces of the facets $\{\simplex_1, \ldots , \simplex_n\}$ become identified and we refer to the equivalence class of such faces as a single face of the triangulation $\tri$. In fact, in $3$-manifold topology it is common to look at so-called {\em $1$-vertex triangulations}, where all $4n$ vertices of its $n$ tetrahedral facets $\{\simplex_1, \ldots , \simplex_n\}$ are identified to a single vertex. The {\em face-vector}, or {\em $f$-vector}, of a PL triangulation $\tri$ is the vector $f(\tri) = (f_0, f_1, \ldots , f_d=n)$ where $f_i$ denotes the number of $i$-dimensional faces of $\tri$.  Given a triangulation $\tri$, its number of top-dimensional simplices $f_d(\tri)$ is referred to as its {\em size} and sometimes also denoted by $n (\tri)$, or just $n$. 

Triangulations of manifolds can also be described in terms of {\it simplicial complexes} (in which any pair of faces can at most have one of their subfaces in common). These are referred to as a {\em combinatorial manifolds} in the literature. Here we mostly focus on the more general notion of PL triangulation defined above, which we denote {\em generalised triangulation}, because it allows us to triangulate a larger variety of manifolds using only a small number of facets. Moreover, given a PL triangulation of a manifold, we can recover a combinatorial manifold of the same PL homeomorphism type by passing to its second derived subdivision.
See \Cref{ssec:Pachner} where we explain how the concept of a PL triangulation naturally leads to a combinatorial description of a PL manifold by its space of triangulations. 

We consider two PL triangulations $\tri$ and $\tri'$ to be equivalent, if one can be transformed into the other by a relabelling of its faces. If this is possible we say $\tri$ and $\tri'$ are {\em (combinatorially) isomorphic}. Given a PL triangulation $\tri$, we can compute its {\em isomorphism signature}. That is, a string of letters encoding the triangulation, such that two triangulations are encoded by the same string if and only if they are isomorphic. Computing the isomorphism signature of an $n$-facet PL triangulation of a $d$-manifold requires $O(d! \cdot n^2)$ time in the worst case \cite{Burton11Simplification}.

\section{MCMC Method: sampling triangulations of manifolds}\label{sec:method}
\label{ssec:MCMC}

Throughout this section, let $\tri$ be a triangulation of a fixed but arbitrary PL, connected, closed $d$-manifold $\manifold$. The goal of our Markov Chain Monte Carlo (MCMC) methods is to sample from the space of triangulations of interest $\Omega$ with a target probability $P(\tri)$. That is, in the limit of large samples $N \rightarrow \infty$ the number of times any $\tri \in \Omega$ is sampled approaches $N_\tri = N \cdot P(\tri)$. Suitable $P(\tri)$'s are chosen to allow for an efficient sample of triangulations of interest. Importantly, the $P(\tri)$'s are known and we can thus also compute estimators based on other $P(\tri)$ -- e.g, a uniform distribution in a subset of $\Omega$ -- by reweighting our samples. This control over the underlying probability $P(\tri)$ of each of the sampled triangulations $\tri$ is essential to control for bias estimates (achieved, e.g., when using naive random-walk explorations of $\Omega$ and counting the number of triangulations with certain properties). 

The MCMC goal of sampling triangulations from $P(\tri)$ is achieved by constructing a Markov chain with transition matrix ${\bf M} = M(\tri \mapsto \tri')$ such that 
\begin{itemize}
    \item[(i)] the probability ${\bf \Pi}$ of each triangulation $\tri$ evolves in time $t$ as 
$${\bf \Pi}_{t+1} = {\bf M} \circ {\bf \Pi}_t.$$
    \item[(ii)] ${\bf \Pi}_t \rightarrow P$ for $t \rightarrow \infty$  for any ${\bf \Pi}_{t=0}$ (e.g., starting at any triangulation $\tri$).
\end{itemize}

In general, sufficient conditions to ensure this convergence is to construct an ${\bf M}$ that~\cite{Robert2004}
\begin{itemize}
\item[(E)] is {\it ergodic} (i.e., irreducible, recurrent, and aperiodic); and
\item[(DB)] satisfies the {\it detailed balance} condition
  \begin{equation}\label{eq.detailedBalance}
    P(\tri) M(\tri \mapsto \tri') = P(\tri') M(\tri' \mapsto \tri) \text{ for all  } \tri,\tri'.
    \end{equation}
  \end{itemize}

The Metropolis-Hastings method~\cite{Robert2004} achieves this goal by decomposing ${\bf M}$ in two steps. In the first step a new triangulation $\tri'$ is proposed. In the second step, it is decided whether this proposed triangulation $\tri'$ is accepted (in which case the state of the chain at $t+1$ becomes $\tri'$) or rejected (in which case the chain remains in $\tri$ at step $t+1$). The transition matrix is thus written as
\begin{equation}\label{eq.MgA}
  M(\tri \mapsto \tri') = g(\tri \mapsto \tri') A(\tri \mapsto \tri'),
  \end{equation}
where $g(\tri \mapsto \tri')$ is the probability of proposing $\tri'$ from $\tri$ and $A(\tri \mapsto \tri')$ is the acceptance probability. 
We can ensure that both \Cref{eq.detailedBalance,eq.MgA} are satisfied by choosing ${\bf A}$ according to the Metropolis choice
\begin{equation}\label{eq.A}
  A(\tri \mapsto \tri') = \min\left(1,\dfrac{P(\tri') g(\tri'\mapsto \tri)}{P(\tri) g(\tri \mapsto \tri')} \right),
\end{equation}
which implies that either $A(\tri \mapsto \tri')=1$ or   $A(\tri \mapsto \tri')=1$ (note that $0 \le {\bf A} \le 1$).

The overall idea of the Metropolis-Hastings MCMC
is that the proposal allows us to perform a random walk in the triangulations $\tri \in \Omega$, which we steer in a controlled way by imposing its convergence to $P(\tri)$ through the acceptance probability~${\bf A}$ in \Cref{eq.A}. Next we discuss suitable choices of proposals $g(\tri \mapsto \tri')$, target probabilities $P(\tri)$, and acceptances~$A(\tri \mapsto \tri')$ -- satisfying relation \Cref{eq.A} -- to the study of triangulations of manifolds.

\subsection{Choice of $P(\tri)$}
\label{ssec:P}

 A primary consideration in our choice of $P(\tri)$ is that the number of  triangulations with $n$ facets (simplices) $|\Omega(n)|$ grows quickly with $n$. We define
\begin{equation}
    |\Omega(n)| \equiv \# \{ \tri \in \Omega \,\mid\, n(\tri) = n \}
\end{equation}
and quantify its growth using
\begin{equation}\label{eq.S}
S^\dagger(n) = \log |\Omega(n)|.
\end{equation}
Since $|\Omega(n)|$ grows quickly with $n$, an unbiased walk very quickly diverges, visiting $\tri$'s with larger and larger $n=n(\tri)$. The choice of $P(\tri)$ and acceptance probabilities has thus one primary task: to counteract the growth of $|\Omega(n)|$ and thus $S^\dagger(n)$ (which acts as an entropic force). 
%  (number of triangulations with the number of facets. 
Our goal is to design $P(\tri)$ to ensure that the random walk generated by the Markov Chain visits triangulations within a characteristic range of sizes of $n$, i.e., with a well defined characteristic size that does not change over a large number of MCMC steps~$t$ when $t \gg 1$.
To achieve this, we choose $P(\tri(n))$ as
\begin{equation}\label{eq.ps}
P(\tri) = \dfrac{e^{-\beta S(n)}}{Z_\beta},
\end{equation}
where $\beta \in \mathbb{R}$ is a parameter, $Z_\beta = \sum_{\tri \in \Omega} e^{\beta S(n)}$ ensures normalisation, % $\sum_{\tri \in \Omega} P(\tri) = 1$, 
and $S(n)$ scales faster than $S^\dagger(n)$ for large $n$. In many situations $|\Omega(n)|$ and thus $S^\dagger(n)$ in \Cref{eq.S} are unknown, but it is enough to have $S(n)$ in \Cref{eq.ps} to grow with $n$ similarly but faster than the growth (of the upper bound) of $S^\dagger$. 

In the case of triangulations, the growth of $\Omega(n)$ is at least exponential (potentially super-exponential) and therefore a safe natural choice that we use in our studies is
\begin{equation}\label{eq.S2}
S(n) = n^2.
\end{equation}
In this case, for sufficiently large $\beta$, the target distribution $P(\tri)$ decays with $n$ faster than $1/|\Omega(n)|$ and the random walk exploration of $\Omega$ does not drift to triangulations $\tri$ with  $n(\tri) \rightarrow \infty $. Intuitively, this imposes a small acceptance probability~(\ref{eq.A}) on bi-stellar moves leading to a larger triangulation, while bi-stellar moves producing a triangulation of the same or smaller size than the current state are always accepted. The move towards larger or smaller $n$ can then be tuned by varying the parameter $\beta$. The constant $Z_\beta$ is irrelevant for the MCMC and does not require computation because the acceptance in \Cref{eq.A} depends only on the ratio of $P(\tri)$ and $P(\tri')$ (i.e., the term $Z_\beta$ cancels out). For instance, considering two triangulations that differ by $\delta$ top-dimensional simplices, $n(\tri') = n(\tri)+\delta$, we obtain from \Cref{eq.S2}
\begin{equation}\label{eq.PratioS}
    \dfrac{P(\tri')}{P(\tri)} = e^{-\beta ((n+\delta)^2-n^2)} = e^{-\beta(2n\delta+\delta^2)},
\end{equation}
which decays exponential with $n$ for fixed $\delta$.

\subsection{Proposal: Bi-stellar moves}% and the Pachner graph}
\label{ssec:Pachner}
The two previously-stated conditions for the convergence of the Metropolis-Hastings MCMC towards $P(\tri)$  --  ergodicity and detailed balance -- impose constraints on the possible proposals $g(\tri\mapsto \tri')$ we can construct. Detailed balance requires that if the proposal admits moves from $\tri$ to $\tri'$ then moves from $\tri'$ to $\tri$ must also be admissible, 
\begin{equation}\label{eq.symmetry}
g(\tri \mapsto \tri')>0 \Leftrightarrow g(\tri' \mapsto \tri)>0,
\end{equation}
Ergodicity translates into a non-zero probability of reaching any $\tri' \in \Omega$  starting from any $\tri \in \Omega$ for $t\rightarrow \infty$. In addition to these restrictions, a computationally efficient proposal must ensure that, given a triangulation $\tri \in \Omega$, all proposed triangulations $\tri'$ are also in $\Omega$. Here we construct the proposal $g(\tri' \mapsto \tri)$ using a well-known construction of how to move through triangulations: given a PL triangulation $\tri$, 
% of a PL $d$-manifold $\manifold$, 
we can pass to a different triangulation $\tri'$ %of %$\manifold$ by the 
using the following local operation: 
\begin{enumerate}
  \item Locate a $(d-i)$-dimensional face $\tau$ surrounded by a set of $(i+1)$ distinct facets. 
  \item Without additional identifications coming from $\tri$ and away from $\tau$, this set of facets forms a subset of the triangulation that is the set of proper faces of a $(d+1)$ simplex. 
  \item Take the complementing $(d+1-i)$ facets of the proper faces of this $(d+1)$-simplex. Naturally, they share a unique common $i$-face. 
  \item Replace the $(i+1)$ facets around the existing $(d-i)$-face of $\tri$ by the $(d+1-i)$ facets around the new $i$-face to form a new triangulation $\tri'$.
\end{enumerate} 

\begin{figure}[htb]
  \includegraphics[width=\textwidth]{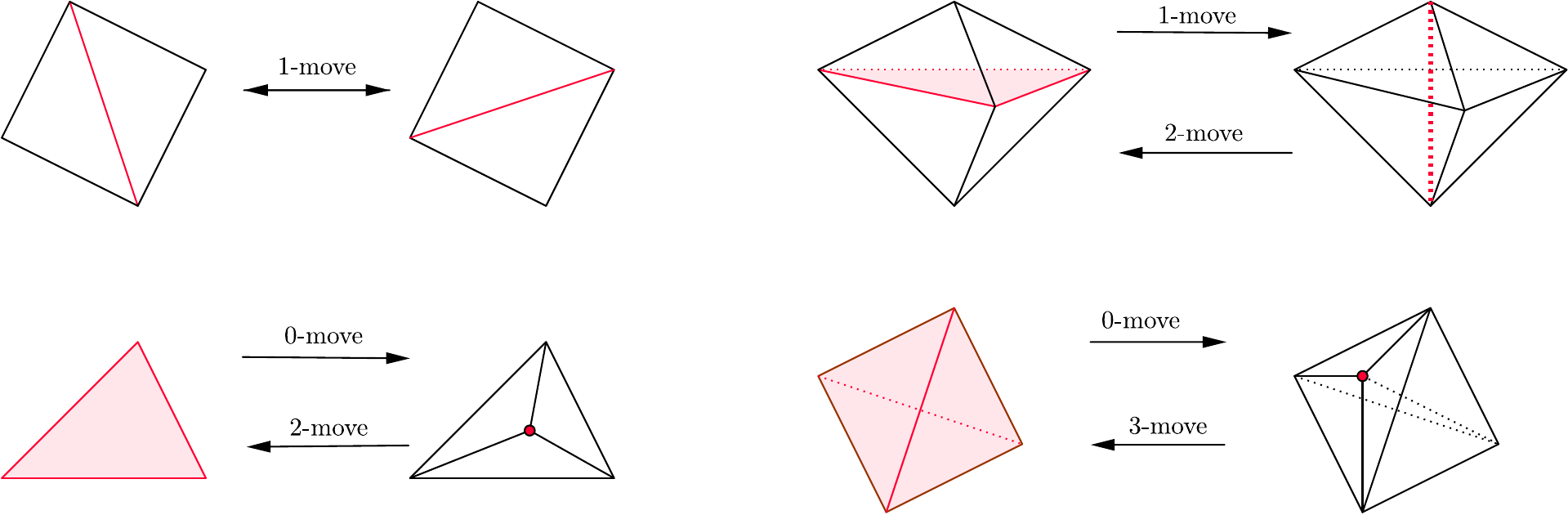}
  \caption{Bi-stellar moves in dimensions two (left) and three (right). For $d=2$, the number $n$ of triangles changes by $(2,0,-2)$ through the application of the $i=(0,1,2)$ move. For $d=3$, the number of tetrahedra $n$ changes by $(3,1,-1,-3)$ through the application of the $i=(0,1,2,3)$ move.
 \label{fig:bistellarmoves}}
\end{figure}

\Cref{fig:bistellarmoves} lists all moves obtained from this local operation in dimensions two and three. It is usually called a {\em bi-stellar $i$-move}, in honour of the dimension of the inserted face -- and accounting for the fact that it describes the stellar subdivision of a $(d-i)$-face followed by the inverse of a stellar subdivision of an $i$-face.  It is sometimes also called a $(i+1)$-$(d+1-i)$-move, indicating how many old facets are replaced by how many new facets. The inverse of an $i$-move exists and is a $(d-i)$-move (and vice versa). Connecting two triangulations if they can be turned into each other by a single move we obtain what is called the {\em Pachner graph}. Our MCMC method performs a random walk in this graph.

It is straightforward to see that, if a PL triangulation $\tri'$ is constructed from a PL triangulation $\tri$ by a bi-stellar move, then $\tri$ and $\tri'$ describe the same PL manifold. But the converse is also true, due to the following seminal theorem by Pachner\footnote{\Cref{thm:Pachner} was initially stated for combinatorial manifolds only, but adapting the result to our more general triangulations is straightforward.}.
\begin{theorem}[Theorem 1 in \cite{Pachner87KonstrMethKombHomeo}]
  \label{thm:Pachner}
  Given two PL triangulations $\tri $ and $\tri'$. Then there exists a sequence of bi-stellar moves turning $\tri$ into $\tri'$ if and only if $\tri$ and $\tri'$ describe the same PL manifold.
\end{theorem}
\Cref{thm:Pachner} ensures the ergodicity of any proposal that assigns a non-zero probability for all bi-stellar moves, i.e. proposals for which all $\tri'$ obtained from bi-stellar moves from $\tri$ have $g(\tri \mapsto \tri') \neq 0$ (the Pachner graph is  connected). Since for any bi-stellar move there is an inverse bi-stellar move, this choice of proposal satisfies \Cref{eq.symmetry} and thus it ensures that detailed-balance is achieved applying the acceptance~(\ref{eq.A}).

In practice, we propose $\tri'$ from $\tri$ following a three step procedure:
\begin{itemize}
    \item[a)] Decide whether $\tri'$ has a larger (up), smaller (down), or the same (stay) $n$ as $\tri$ (the last option is only possible in even dimensions by a straightforward parity argument, see also \Cref{fig:bistellarmoves}).
    \item[b)] With uniform probability, choose one of the types $i$ of bi-stellar moves that satisfies the decision in a) (up/down/stay).
    \item[c)] With uniform probability, pick one $\tri'$ obtained from $\tri$ by performing an $i$-move (i.e., same probability for all $i$-neighbours of $\tri$ in the Pachner graph, see discussion below). 
\end{itemize}
Step a) provides another opportunity to favour triangulations of interest (e.g., of smaller $n(\tri)$), a decision that impacts the mixing time and convergence of the Markov Chain to $P(\tri)$. We typically choose the probability $\alpha$ of going up to be such that $\alpha \rightarrow 0$ when $n\rightarrow \infty$.  This takes advantage of the fact that bi-stellar moves are local in $n$ -- they change $n$ at most by $\pm d$ (typically $d\ll n$) -- and that we know in advance how each move affects $n$. In step b), recall that performing a bi-stellar $i$-move on $\tri$ produces a triangulation $\tri'$ with more facets for $i<d/2$, with fewer facets for $i > d/2$, and with the same number of facets as $\tri$ for $i=d/2$.  
Step c) requires the enumeration of all isomorphism types of triangulations that result from performing a bi-stellar $i$-move on $\tri$. For every $0 \leq i \leq d$, this produces a list of triangulations that we call the {\em $i$-neighbours of $\tri$}. Here, we specifically exclude $\tri$ from the set of $i$-neighbours.
Note that $i$-neighbours and $j$-neighbours have distinct $f$-vectors for $i\neq j$, and hence every combinatorial isomorphism type of triangulation of $\manifold$ can only occur as one type of neighbour of $\tri$. By the definition of a bi-stellar $i$-move, the number of $i$-neighbours is trivially bounded from above by $f_{d-i}(\tri)$. Moreover, the number of $i$-neighbours is typically very close to $f_{d-i}(\tri)$ for $i \approx 0$, and very small for $i \approx d$ (and in between for $0 \ll i \ll d$).
Therefore, one safe procedure to ensure that $i$-neighbours of $\tri$ are chosen with uniform probability is to attribute a probability given by $1/f_{d-i}(\tri)$ for $i<d/2$ or $1/f_{i}(\tri')$ for $i\ge d/2$, 
because $f_{d-i}$ (respectively, $f_i$) is the maximum number of possible faces in which we can apply the move $i$.  Typically this leaves some unallocated probability (i.e., no independent $\tri'$ to be chosen), in which case we choose $\tri'=\tri$. 

We now compute the proposal probability $g(\tri \mapsto \tri')$ obtained from the combination of steps a), b), and c) described above. Note that the choices in each step are independent and that a given $\tri'$ can only be reached from $\tri$ through one of the $i$ moves.  The proposal probability $g(\tri \mapsto \tri')$ is thus obtained as the product of the steps a), b), and c) as
%
%\begin{equation}\label{eq.gi} 
\begin{eqnarray}\label{eq.gi} 
g (\tri\mapsto \tri') = \left. \begin{cases}  \dfrac{\alpha(n)}{\lfloor d/2\rfloor f_{d-i} (\tri)} & \text{for } i < \dfrac{d}{2} \Rightarrow n'>n \text{ (move up in $n$),}\\ 
& \\
\dfrac{(1-\alpha(n)-\tilde{\alpha}(n))}{f_{d/2} (\tri)} & \text{for } i = \dfrac{d}{2} \Rightarrow  n'=n \text{ ($d$ even only),} \\  
& \\
\dfrac{\tilde{\alpha}(n)}{\lfloor d/2\rfloor f_{i} (\tri')} & \text{for } i > \dfrac{d}{2} \Rightarrow n'<n  \text{ (move down in $n$),}  \end{cases} \right. \end{eqnarray} 
where $i= 0, 1, \ldots, d$ corresponds to the $i-$th bistellar move for which $\tri'$ is an $i$-neighbour of $\tri$, $f_i$ introduced in \Cref{ssec:triangulations} is the number of $i$-dimensional faces of $\tri$, $0<\alpha(n) <1$ is the probability of choosing an $i$-move with $n'>n$, $0<\tilde{\alpha}(n) <1$ is the probability of choosing an $i$-move with $n'<n$, 
 and $\lfloor \ldots \rfloor$ is the floor function so that $\lfloor d/2 \rfloor$ is the number of $i$-move types that lead to $n'>n$ (or equivalently, $n'<n$). We typically choose 
\begin{equation}\label{eq.alphatilde}
\tilde{\alpha} = (1-\alpha)/r
\end{equation}
where $r$ controls the fraction of non-up moves that go down (for odd $d$ we must have $r=1$, for even $d$, $r > 1$, to ensure a non-zero probability goes to $d/2$-moves). As we see below, a good choice of $\alpha(n)$ to control for the growth of $|\Omega(n)|$ as in \Cref{eq.S2} is to choose
\begin{equation}\label{eq.gamma}
\alpha(n) = e^{-\gamma n},
\end{equation}
for varying $\gamma>0$.

\subsection{Choice of $A(\tri\mapsto\tri')$}
\label{ssec:acceptance}

The probability of accepting a move in the Metropolis-Hastings MCMC is fixed once the target probability $P(\tri)$ and the proposal $g(\tri\mapsto \tri')$ are specified. For instance, using the specific choices above for $P(\tri)$ -- \Cref{eq.S2} and thus \Cref{eq.PratioS} -- and the proposal -- \Cref{eq.alphatilde,eq.gamma} --  we obtain an exact expression for the acceptance in \Cref{eq.A} as computed in \Cref{app.acceptance}. Most importantly, it scales with $n$ as

\begin{equation}\label{eq.Ascaling}
A_\pm \sim e^{\pm(\gamma-2\beta |\delta|) n}, 
\end{equation}
where $\delta = n'-n$ and the "$+$" ("$-$") case corresponds to moving up with $\delta>0$ (down with $\delta<0$).
This suggests that the only choice for which one of the acceptances does not vanish for $n\rightarrow \infty$ is
\begin{equation}\label{eq.gammabeta}
  \gamma = 2\beta |\delta|,
\end{equation}
where $|\delta| \le d$. Choosing a non-vanishing acceptance is crucial for the success of our methods because small acceptance rates increase the correlation between samples and thus the time needed for the chain to relax to the equilibrium distribution $P(\tri)$, removing the advantages of MCMC. 

The reasoning above suggests an alternative approach for setting the MCMC. Instead of fixing $P(\tri)$ we fix the acceptance to the ideal case $A=1$, maintaining $g(\tri \mapsto \tri')$ as in \Cref{eq.gi} with \Cref{eq.alphatilde} and \Cref{eq.gamma}. We can then compute the relative change of $P(\tri)$ for this case inverting \Cref{eq.A} as 
\footnote{This is obtained introducing \Cref{eq.gratioU} (the ratio of proposals $g$ for moving up) into \Cref{eq.A}.  An equivalent result is obtained considering the case of moving down, introducing \Cref{eq.gratioD} in \Cref{eq.A} to obtain
\begin{equation}\label{eq.Adown}
    \dfrac{P(n+\delta)}{P(n)} = \dfrac{\tilde{\alpha}(n)}{\alpha(n+\delta)} = \frac{1}{r} \dfrac{1-\alpha(n)}{\alpha(n+\delta)},
\end{equation}
which recovers the previous case considering $n=n'+\delta'$ with $\delta'=-\delta$.}
\begin{eqnarray}
    \dfrac{P(\tri')}{P(\tri)} =  & \dfrac{P(n+\delta)}{P(n)}  = & \dfrac{g(\tri \mapsto \tri')}{g(\tri' \mapsto \tri)}=  \dfrac{\alpha(n)}{\tilde{\alpha}(n+\delta)} = \nonumber \\ 
    = & r \dfrac{\alpha(n)}{1-\alpha(n+\delta)}= & r \dfrac{e^{-\gamma n}}{1-e^{-\gamma(n+\delta)}} (\approx re^{-\gamma n} \text{ for } n \gg 1/\gamma)
    \label{eqPratioS2}
\end{eqnarray}
For large $n$, we recover the exponential decay in $n$ obtained in \Cref{eq.PratioS} and, comparing the two expressions, we retrieve the relationship between $\gamma$ and $\beta$ obtained in \Cref{eq.gammabeta}. The advantage in this case is that $\gamma$ is the only control parameter (which controls the growth in $|\Omega(n)|$) and that all moves are accepted, leading to a faster convergence of the MCMC. A potential drawback of this {\it accept all} method is that there is no explicit expression for $P(n(\tri))$, in contrast to a method based on the choice of $\beta$ in \Cref{eq.ps}. Instead it needs to be computed iteratively from Equation~(\ref{eqPratioS2}) (e.g., fixing $P(\delta)$ and imposing normalisation $\sum_n P(n)=1$). However, as argued before in the computation of $Z_\beta$, the exact probability $P(\tri)$ is typically not achieved and the ratio of $P$s (relative probability) is typically sufficient for reweighting the numerical results and estimating quantities of interest. 

\subsection{Implementation}

{\it In principle}, the method discussed above is applicable to arbitrary manifolds and dimensions~$d$. {\it In practice},  adaptations are necessary to address challenges and opportunities that are specific to each dimension $d$. 
Below we describe the implementation in $d=2$ and $d=3$, and discuss choices and challenges that appear for $d\ge4$. Details can be found in the \Cref{app.numerics} and the Python code in the repository~\cite{gitHub}.

\subsubsection{Dimension 2} \label{ssec:dim2Desc}
Throughout this section, let $\tri$ be a triangulation of a connected, closed surface $\surface$. We denote the Euler characteristic of $\surface$ by $\chi = \chi (\surface)$. The $f$-vector of a triangulation of $\surface$ is 
$$f(\tri) = (f_0,f_1,f_2)=(n/2+\chi, 3n/2, n).$$ 
Note that this implies that the number of triangles $n$ of a triangulation of a surface must be a positive even number.
It also implies the sharp lower bound $n \geq \max \{2, 2 (1-\chi)\}$. 

The Pachner graph of surfaces has been extensively studied in the literature. Fixing $\surface$, every pair of (abstract) triangulations $\tri$ and $\tri'$ of $\surface$ can be connected by a sequence of bi-stellar moves of length $O(n(\tri) + n(\tri'))$ -- where the constant only depends on the genus $g$ of $\surface$.\footnote{Finding such a sequence is straightforward: Both $\tri$ and $\tri'$ must have a vertex contained in only $O(g)$ triangles. Such a vertex can then be removed using only $O(g)$ bi-stellar moves. Since there are $O(n(\tri))$ ($O(n(\tri'))$) vertices in $\tri$ ($\tri'$), iterating this procedure on both $\tri$ and $\tri'$ produces two $1$-vertex triangulations of $\surface$ in  $O(g(n(\tri) + n(\tri')))$ steps. These $1$-vertex triangulations -- each with $O(g)$ triangles -- are connected by a sequence of $1$-moves of length a function of $g$.} Moreover, if $\tri$ and $\tri'$ have the same number of triangles, then $\tri$ can always be transformed into $\tri'$ by a sequence of $1$-moves. In other words, the flip graph of $n$-triangle triangulations of $\surface$ is connected.
If $\surface$ is a $2$-dimensional sphere, such a sequence of length at most $5 \cdot n + O(1)$ always exists~\cite{CARDINAL2018206}. For triangulations of surfaces of arbitrary genus $g$, the maximal length of a sequence of $1$-moves between any two $n$-triangle triangulations is at most $O(g \log(g+1) + n \log (n))$, as can be followed from \cite[Theorem 1.4]{Disarlo14FlipDistances} -- a bound that even applies to the more challenging setting where vertices have a fixed labelling. These facts motivate the implementation described in \Cref{algo:dim2} for $d=2$.

\begin{algorithm}[htb]
\caption{Metropolis-Hastings MCMC to sample triangulations of manifolds in $d=2$, a particular case of the {\it accept all} method described in \Cref{ssec:acceptance}.}\label{algo:dim2}
\begin{description}
  \item[0.] Input:\\  
    Triangulation $\tri$ of a surface $\surface$, \\
    $\gamma >0$, we use $ \gamma= 1/k$ for increasing $k \in \mathbb{N}$,  \\
    t=0
  \item[1.] Sample $u \in U([0,1])$:
    \\ If $u < e^{-\gamma n(\tri)}$:\qquad\qquad\qquad\qquad\   $i:=0$, $m:=n(\tri)$
    \\If $e^{-\gamma n(\tri)} \leq u \leq \frac{1+e^{-\gamma n(\tri)}}{2}$:    \qquad\ $i:=1$, $m:=\frac{3}{2}n(\tri)$
    \\ If $u \geq \frac{1+e^{-\gamma n(\tri)}}{2}$:\qquad\qquad\qquad\quad\      $i:=2$, $m := n(\tri)-2$
  \item[2.] Enumerate $i$-neighbours $\tri_1, \ldots , \tri_\ell$ of $\tri$
  \item[3.] Sample $v \in U([0,1])$:
\\If $v < \ell / m$:\quad pick $\tri'$ from $\tri_1, \ldots , \tri_\ell$ at random
\\If $v \geq \ell/m$:\quad pick $\tri' = \tri$
  \item[4.]  Sample $\tri'$, update state ($\tri := \tri'$, $t:=t+1$), and go to {\bf 1.}
\end{description}

\end{algorithm}

\subsubsection{Dimension 3} 
\label{ssec:dim3Desc}

\begin{algorithm}[htb]
\caption{Metropolis-Hastings MCMC to sample triangulations of manifolds in $d=3$, a particular case of the {\it accept all} method described in \Cref{ssec:acceptance}. \label{algo:dim3}}
\begin{description}
  \item[0.] Input: 
  \\Triangulation $\tri$ of a $3$-manifold $\manifold$
  \\$\gamma = 1/k$, $k \in \mathbb{N}$  
  \\t=0
  \item[1.] Sample $u \in U([0,1])$
  \\If $u < e^{-\gamma n(\tri)}$:\qquad\qquad\ $i:=1$, $m:=2n$
  \\ If $u \geq e^{-\gamma n(\tri)}$:\qquad\qquad\ $i:=2$, $m:=2n-2$
  \item[2.] Enumerate $i$-neighbours $\tri_1, \ldots , \tri_\ell$ of $\tri$
  \item[3.] Sample $v \in U([0,1])$
  \\ If $v < \ell / m$:\quad pick candidate $\tri'$ from $\tri_1, \ldots , \tri_\ell$ at random
  \\ If $v \geq \ell/m$:\quad set $\tri' = \tri$
  \item[4.]  Sample $\tri'$, update state ($\tri := \tri'$, $t:=t+1$), and go to {\bf 1.}
\end{description}
\end{algorithm}

Throughout this section, let $\tri$ be a triangulation of some fixed connected, closed $3$-dimensional manifold $\manifold$. The face vector is given by
\begin{equation}
\label{eq:dim3}
f(\tri) = (f_0,f_1,f_2,f_3=n)=(f_0,n+f_0,2n,n).
  \end{equation}
Much less is know about the Pachner graphs obtained in these cases. For instance, it is not known whether the number of triangulations of the $3$-dimensional sphere is singly- or super-exponential in the number of tetrahedra. Given a pair of triangulations of a $3$-manifold, upper bounds for the lengths of shortest sequences of bi-stellar moves between them exist, but are believed to be far from optimal \cite{Mijatovic03BoundsS3,Mijatovic04SFS}. While some of the bounds in \cite{Mijatovic04SFS} are towers of exponentials, it is believed that sharp upper bounds exist for, say, the $3$-dimensional sphere that are polynomial in the number of tetrahedra of the input triangulation: in the related problem of untangling diagrams of the unknot to the trivial diagram, a number of Reidemeister moves polynomial in the initial crossing number of the diagram is sufficient \cite{Lackenby15ReidemeisterMoves}.

There are, however, two results in the literature that provide a pathway to a tractable method in the $3$-dimensional setting.
\begin{enumerate}
  \item There exists a procedure to transform an arbitrary $n$-tetrahedra triangulation of a $3$-manifold (satisfying some weak topological assumptions) into one that only uses one vertex, and no more than $n$ tetrahedra. Moreover, this procedure uses a number of local moves polynomial in $n$ \cite{burton14-crushing-dcg,burton12-unknot}.
  \item Any two one-vertex triangulations $\tri$ and $\tri'$ of $\manifold$ (each with at least two tetrahedra) are connected by a sequence of only $1$- and $2$-moves \cite{MatveevSpecialSpines}.
\end{enumerate}

In order to turn the task of searching through the space of triangulations of a $3$-manifold into a feasible problem, we make use of these two facts and focus on the space of $1$-vertex triangulations (with face-vectors $(1,n+1,2n,n)$), connected only by $1$- and $2$-moves. More precisely, we propose the method given in \Cref{algo:dim3}. From an optimisation standpoint, our restriction to $1$-vertex triangulations is not severe because, with very few known exceptions, the smallest triangulations of a $3$-manifold always have only one vertex \cite{burton14-crushing-dcg,burton12-unknot}.

\subsubsection{Dimension $\ge 4$}
\label{ssec:dim4Desc}

Exploring the space of triangulations of a $d$-manifold, $d\geq 4$, is much harder than in lower dimensions: some pairs of small triangulations of $4$-spheres require connecting sequences of moves to pass through considerably larger intermediate triangulations. Such a phenomenon does not exist in dimension two, where all $n$-triangle triangulations of a surface are connected via $1$-moves, and it is not known to exist in dimension three \cite{Burton11Simplification}.
This behaviour is not surprising, given that the homeomorphism problem -- the task of deciding whether two given (triangulations of) $d$-manifolds are topologically equivalent or not -- becomes undecidable for $d\geq 4$ \cite{Markov58HomeoProb}.

These points indicate that our MCMC method in $d\ge4$, while straightforward to implement, is expected to take longer to converge and to show larger correlations between samples. One approach to mitigate such negative effects (potentially already at $d=3$) is to use different types of proposals, beyond just bi-stellar moves.
Another possibility is to restrict the random walk to triangulations with a fixed number of vertices (as done in \Cref{algo:dim3}). Such a procedure only performs $1$-, $2$-, and $3$-moves and conceptually resembles \Cref{algo:dim2}. But the restriction to a fixed number of vertices is much more severe in dimension four than it is in dimension three: we loose  guarantees for ergodicity, and, for instance, the smallest possible triangulation of a given manifold does not necessarily have only one vertex (on the contrary, in many cases $1$-vertex triangulations are ``unnecessarily'' large \cite{Tobin234Manifolds}). 

\section{Numerical experiments}
\label{sec:experiments}

In this section we report on numerical experiments obtained using our MCMC method in dimensions 2 and 3. We start by showing that the (asymptotic) theoretical properties of the method are observed in finite simulations. We then reproduce exact known results in dimension 2 to confirm the accuracy of the method. Finally, we report on new findings in dimensions 2 and 3. In all applications we use our accept-all method described in \Cref{ssec:acceptance}, with parameters $\alpha$ and $\tilde{\alpha}$ given by \Cref{eq.alphatilde,eq.gamma}, $r=1$ (for odd $d$) or $r=2$ (for even $d$), and varying $\gamma$. We typically choose $\gamma= \frac{1}{k}$, $k \in \{1, 2, 3, \ldots 25\}$, and $20$ separate runs of $10\ 000\ 000$ moves ($100\ 000$ samples each, that is, collecting every $100$th triangulation in the MCMC chain).
%%%%%%%%%%%%%%%%%%%%%%%%%%%%%%% Numerical experiments %%%%%%%%%%%%%%%%%%%%%%%%%%

\begin{figure}[hbt]
  \centerline{
  \includegraphics[height=6cm]{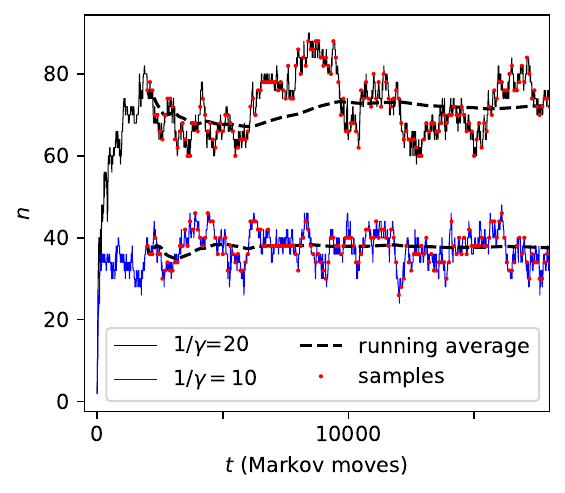}\hspace{-0.6cm}\includegraphics[height=6cm]{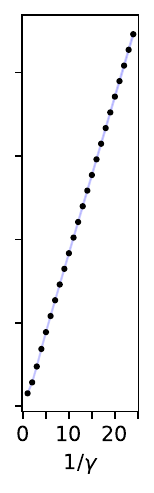}\includegraphics[height=6cm]{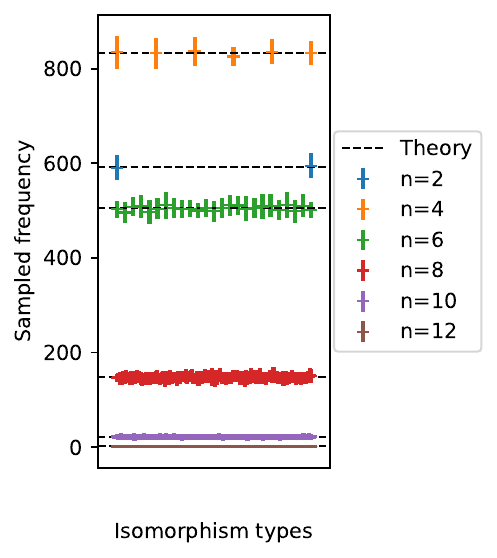}
  }
    \caption{Our MCMC method samples triangulations of various sizes $n$ uniformly. Left: Convergence to characteristic triangulation sizes $n$ as a 
    function of Markov time $t$ for two values of the parameter $\gamma \in \{ \frac{1}{10}, \frac{1}{20}\}$ in Algorithm~\ref{algo:dim2}
    (see legend). At $t=0$, we start at a $n=2$-triangle seed triangulation of the $2$-sphere obtained by gluing two triangles along their boundaries. 
    We estimate the average size $\langle n \rangle$ sampling every $100$ iterations for $t\in[2\ 000, 100\ 000]$. 
    Centre: $\langle n \rangle$ as a function of $1/\gamma$.
    Right: triangulations with $n$ triangles are sampled with the same probability. The sampled frequencies (y-axis) of individual isomorphism types 
    of triangulations (ordered along the x-axis) is shown for different sizes of $n$ (see legend).
    Numerical results were obtained using  $\gamma=\frac13$ and $\approx 90000$ samples in each realisation (average and standard deviation over $20$ 
    realisations are reported). The theoretical expectations (dashed lines) are the number of samples of size $n$ divided by the (known) number of 
    different $n$-triangle triangulations. 
    %on average $1186$, $4998$, $12626$, $23000$, $27110$, $19919$ triangulations of size $n = 2,4,6,8,10,12$, respectively.  
 }
    \label{fig:firstNumerics}
\end{figure}

\subsection{Properties of the method}
\label{ssec:basic}

Our first numerical experiments aim to confirm that our MCMC achieves in practice what it was designed to achieve in \Cref{sec:method}. 
One critical aspect is the ability to control for the growth in the number of triangles $n$ seen in the chain. The numerical results 
shown in \Cref{fig:firstNumerics}(left and centre panels) confirm that our MCMC method succeeds in this task: for a given value of the MCMC control parameter $\gamma$, the size $n$  of the triangulations eventually stops growing and starts oscillating around a typical value of $n$; and by changing $\gamma$ we can choose the $n$ that we predominantly sample.  We also test that our method samples uniformly triangulations with the same $n$.  For small $n$, the complete list of all isomorphism types of $n$-triangle triangulations of the $2$-sphere is known and can thus be used to test our method. The results shown in \Cref{fig:firstNumerics} on the right confirm that our method succeeds not only in finding these triangulations but also in sampling them with frequencies that depend only on $n(\tri)$ -- according to $P(n)$ -- so that different triangulations $\tri$ with the same $n$ are sampled with the same frequency.

\subsection{Number of triangulations}
\label{ssec:numberoftriangulations}

We now study the growth of the number of triangulations $|\Omega(n)|$ with triangulation size $n$. Our focus is on determining whether an exponential scaling $|\Omega(n)|\sim e^{\tilde{c} n}$ holds for $n\rightarrow \infty$, in which case the following ratio converges to a constant
\begin{equation}\label{eq.R}
R\equiv \frac{|\Omega(n+\tilde{\delta})|}{|\Omega(n)|} \quad \rightarrow \quad \, C\,= e^{\tilde{\delta} \tilde{c}},
\end{equation}
where $\tilde{\delta}=2$ ($\tilde{\delta}=1$) for even (odd) dimensions.  For every choice of the MCMC parameter $\gamma$ we can precisely compute the bias of the MCMC in choosing a triangulation of size $n$ over a triangulation of size $n + \tilde{\delta}$, which is then accounted for to obtain an unbiased estimator $\hat{R}$ of $R$ (see \Cref{app.reweighting} for details of this reweighting process). The absolute number $|\Omega(n)|$ can be estimated from $\hat{R}$ using the known values of $|\Omega(n)|$ for small $n$, e.g., $|\Omega(4)|=1$ for simplicial triangulations of the $2$-sphere.

\subsubsection{The number of simplicial triangulations of the $2$-sphere} 
\label{ssec:simplicialS2}

We start testing whether our method can numerically reproduce a theoretical result about the number of {\it simplicial} triangulations of the $2$-sphere, essentially due to Tutte \cite{tutte_1962}. In his work, Tutte proved that the number of rooted, $n$-triangle simplicial triangulations of the $2$-sphere, up to combinatorial isomorphism, asymptotically equals
$ \frac{3}{16\sqrt{6 \pi \left (\frac{n}{2} + 1\right )^5}} \left ( \frac{256}{27} \right )^{n/2}$.
Here, {\em rooted} means that a single triangle of the $2$-sphere together with an ordering of its vertices is taken as a frame of reference against which triangulations are compared for combinatorial isomorphy. Taking into account that for each $n$-triangle triangulation of the $2$-sphere we can do this in $3! \cdot n = 6n$ ways, this means that the overall number of $2$-sphere triangulations up to isomorphism must, asymptotically, be bounded from below by 
\begin{equation}
  \label{eq:tutte}
  |\Omega(n)| \ge \frac{1}{6n}\cdot \frac{3}{16\sqrt{6 \pi \left (\frac{n}{2} + 1\right )^5}} \left ( \frac{256}{27} \right )^{n/2} .
\end{equation}
For large $n$, all $6n$ representations of a triangulation become distinct (because symmetries in large triangulations are very rare) and hence this bound converges to the actual number of isomorphism types of simplicial triangulations. Regardless of these polynomial pre-factors, the ratio $R$ from \Cref{eq.R} converges to $C=256/27$.

\begin{figure}[htbp]
    \centerline{\includegraphics[width=1\textwidth]{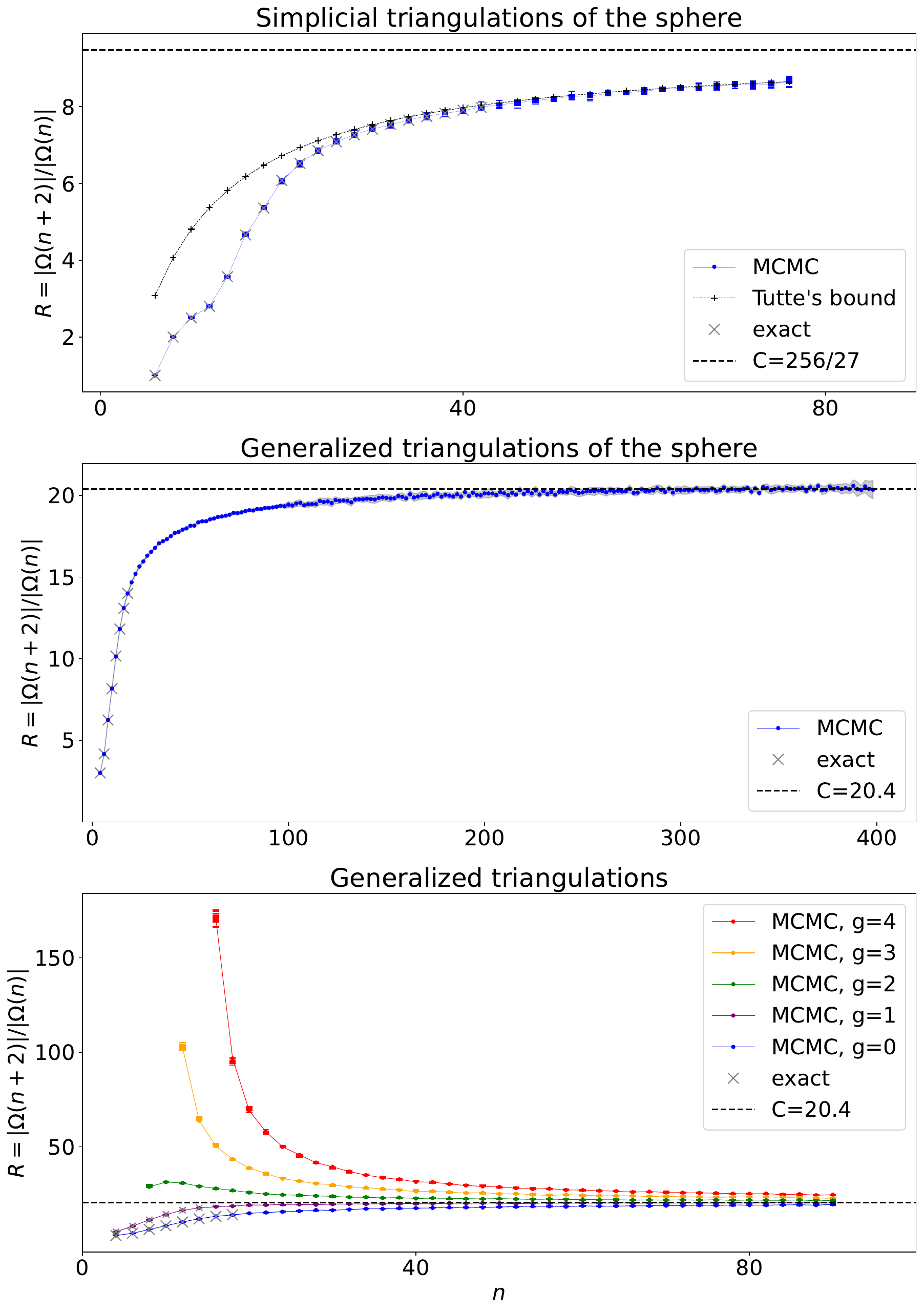}}
  \caption{Rate of growth of the number of isomorphism types of triangulations $|\Omega(n)|$ for manifolds in $d=2$ dimensions and different genera $g$. The y-axis shows $R=|\Omega(n+2)|/|\Omega(n)|$ and the x-axis the size of the triangulation $n$ (number of triangles). 
    Our MCMC results agree (within $99\%$ confidence intervals) with the exact results obtained from enumeration for small $n$~\cite{ManifoldPage,Burton11Simplification}, agree with the sharp asymptotic bound due to Tutte in \Cref{eq:tutte} in the simplicial case for $d=2$ spheres, and show a convergence to a constant $C$ as indicated in \Cref{eq.R}. 
   In all cases we use $\gamma = 1/k, k \in \{1, 2, \ldots, 25\}$ and for the $2$-sphere (middle panel) additionally $k \in \{0,35, \ldots , 100\}$.  Symbols show the average and variation over $20$ separate runs of $10^7$ moves ($10^5$ samples) each.}\label{fig:R}
\end{figure}

In order to experimentally reproduce this result, we ensure that the space $\Omega$ explored by our MCMC corresponds to the space of {\it simplicial} triangulations by using a slight variation of \Cref{algo:dim2}: we start with a simplicial triangulation of the $2$-sphere, and we only allow bi-stellar moves producing a neighbouring simplicial triangulation.\footnote{This is the original setting Pachner worked with when proving \Cref{thm:Pachner}, \cite{Pachner87KonstrMethKombHomeo}} 
The results we obtain are summarised in \Cref{fig:R}(top panel) and show excellent agreement not only with Tutte's asymptotic formula from \Cref{eq:tutte}, but also with exact data from enumeration for small $n$.

\begin{figure}[htb]
    \centerline{\includegraphics[width=1\textwidth]{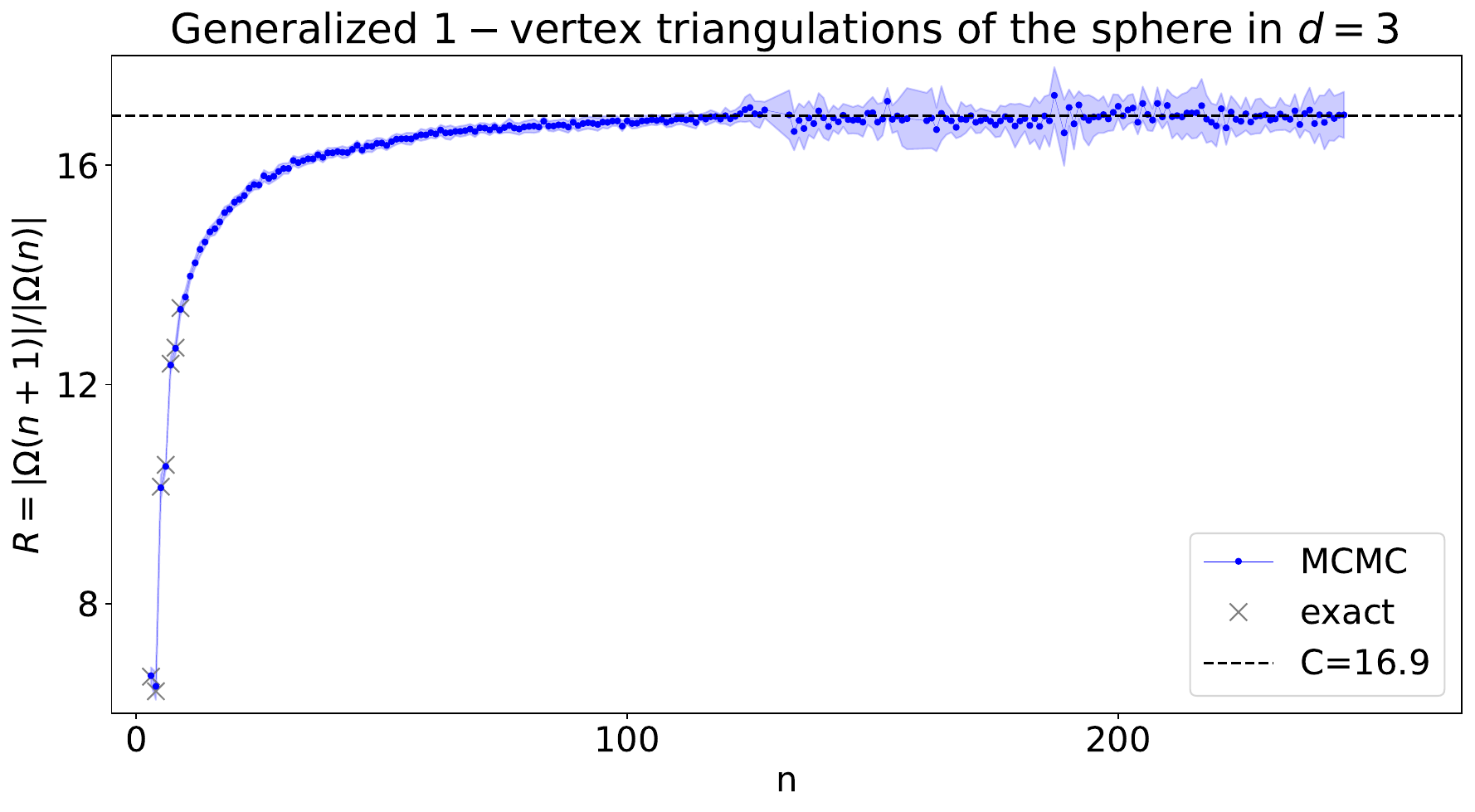}}
  \caption{Rate of growth of the number of isomorphism types of triangulations $|\Omega(n)|$ for $d=3$ dimensional spheres. The y-axis shows $R=|\Omega(n+1)|/|\Omega(n)|$ and the x-axis the size of the triangulation $n$ (number of tetrahedra). 
    Our MCMC results agree (within $99\%$ confidence intervals) with the exact results obtained from enumeration for small $n$~\cite{Burton11Simplification} and show a convergence to a constant $C$ as indicated in \Cref{eq.R}. 
   We use $\gamma = 1/k, k \in \{1, 2, \ldots, 40,50, \ldots ,90,100\}$. Symbols show the average and confidence interval (standard deviation of the mean) obtained over $20$ separate runs of $10^7$ moves ($10^5$ samples) each.}\label{fig:R3}
\end{figure}

\subsubsection{Generalised triangulations of higher-genus surfaces}
\label{ssec:higherGenus}

We next investigate generalised triangulations of surfaces $\surface$ of genus $0\leq g \leq 4$. Let $\Omega(g,n)$ denote the set of isomorphism classes of $n$-triangle generalised triangulations of the closed, orientable surface of genus $g$. 
Recall that a triangulation of a closed, orientable surface $g$, $g>0$, requires at least $n\geq 4g-2$ triangles. 
It is known that, for $g$ fixed and {\em simplicial} triangulations of surfaces, 
the number of triangulations grows exponentially with $n$ so that $|\Omega(g,n+2)| / |\Omega(g,n)| \to C$ for some constant $C \in \mathbb{R}$ not depending on $g$ 
\cite{Gao91NumberOfRootedTriangulations,Goulden08KPHierarchy,Bender13ApproximationNumberTriangulations}.
In this section we numerically investigate this behaviour in the case of {\it generalised} triangulations. 

Our numerical results are obtained this time allowing for all bi-stellar moves and thus sampling from 
the space $\Omega$ of generalised triangulations $\tri$. Results summarised in \Cref{fig:R}(middle panel) show for the $2$-sphere ($g=0$) that our MCMC estimates agree with enumeration data (for small $n$) and allow us to obtain estimates for 
larger triangulations, which show that the growth stops accelerating and seems to converge to a constant. 
We estimate from this plot $C=20.43\pm0.05$ for this constant,  which implies from \Cref{eq.R} an exponential growth of $|\Omega(n)|$ with a rate $\tilde{c} = 1.5085 \pm 0.0012$.
Note that this $C$ is substantially larger than the one for {\it simplicial} triangulations in the previous experiment, 
with the difference providing a quantification of how scarce simplicial triangulations are in the set of all generalised 
triangulations of the $2$-sphere. As we can observe from \Cref{fig:R} (lower panel), the differences in growth rates of numbers of triangulations of surfaces of genus $0$ to $4$ differ 
significantly for small values of $n$, but then all converge to the same constant $C \approx 20.4$ until they cannot be distinguished from 
each other anymore numerically. This is in line with research presented in \cite{Gao91NumberOfRootedTriangulations,Goulden08KPHierarchy,Bender13ApproximationNumberTriangulations}.

\subsubsection{The number of triangulations of the 3-sphere}
\label{ssec:dim3}

Much less is known about the growth of isomorphism types of $n$-tetrahedra triangulations of a fixed $3$-manifold $\manifold$. The question of whether the number of isomorphism types of $n$-tetrahedra triangulations of the $3$-sphere is singly exponential, or super-exponential in $n$, has been the focus of considerable research efforts. We know that the number of isomorphism types of $n$-tetrahedra triangulations of the $3$-sphere grows at least singly exponential, see, for instance, \cite{Pfeifle02ManySpheres}, and at most at the rate of $O( D^{n \log n})$, for some constant $D \in \mathbb{R}$ \cite{Stanley75UBT}. This, and the success of our method in $d=2$, motivates us to study these questions experimentally.

We focus on the class of generalised triangulations of the $3$-sphere with one vertex. Naturally, a super-exponential growth for the number of $1$-vertex triangulations implies the same for the case of general triangulations. Conversely, it seems unlikely conceptually, that a singly exponential growth rate for $1$-vertex triangulations allows for super-exponential growth in the general setting. 
The results shown in \Cref{fig:R3} show numerical evidence for a single-exponential growth of $|\Omega(n)|$ in $1$-vertex $3$-sphere triangulations: the growth of $R(n)$ shows the same characteristics than the ones we observe in the two cases in dimension two (in which $R\to C$ is known theoretically), and the results we obtain are compatible with a constant $R$ (within our numerical precision, no linear trend in the $R$ vs. $n$ curve is observed for $n\ge 180$).  We estimate from this plot $C =16.89\pm0.05$, which implies from \Cref{eq.R} an exponential growth of $|\Omega(n)|$ with a rate $\tilde{c} = 2.827 \pm 0.003$. While we cannot of course discard the appearance of slow (e.g., $O(\log (n))$) growth in $R$, we can expect this to be visible only for $n\gg 200$.

As a caveat, the ability of our method to efficiently obtain independent samples of triangulations the $3$-sphere is affected by the fact that -- sometimes -- multiple $1$-moves are necessary to connect a given $3$-sphere triangulation to some other triangulations with the same number of tetrahedra. While, in practice, this number of $1$-moves seems to be very small in a typical scenario \cite{Burton11Simplification}, this phenomenon may have an impact on the speed and type of convergence to a uniform sampling procedure. Here, the agreement of estimations based on our method with exact enumerations (for small $n$) provides further evidence for the accuracy of the estimation of confidence intervals based on our methods. 

\begin{figure}[htbp]
  \centerline{\includegraphics[width=0.55\textwidth]{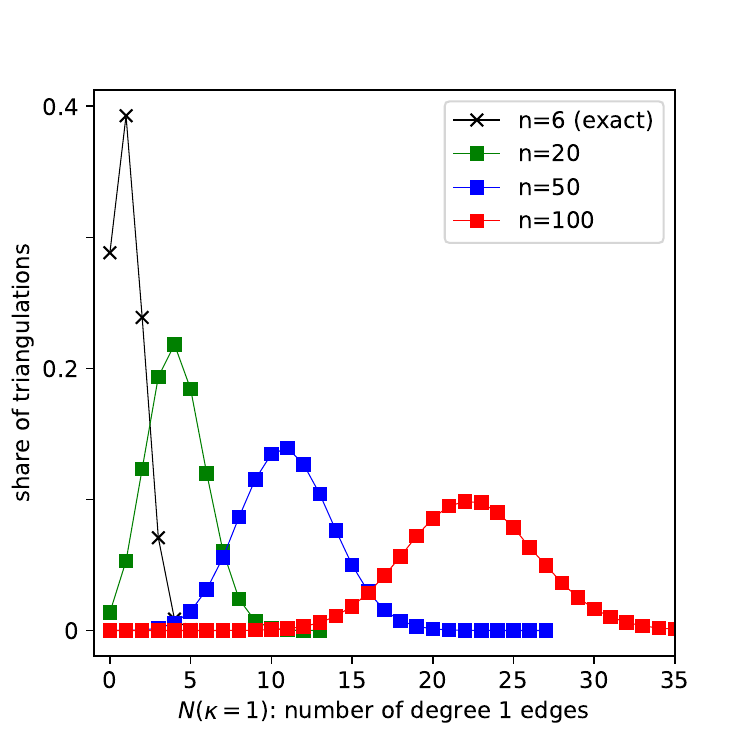}\includegraphics[width=0.55\textwidth]{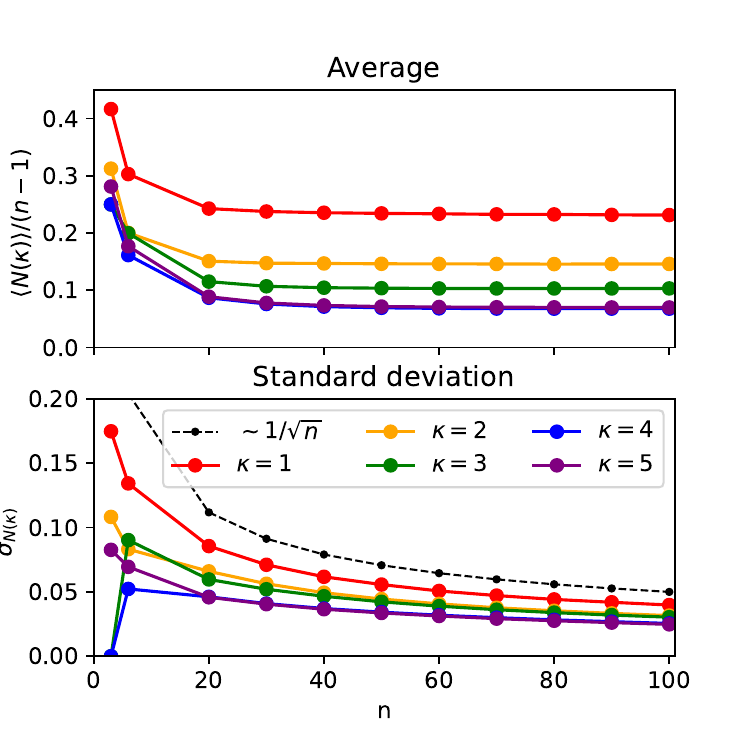}}
  \centerline{\includegraphics[width=0.55\textwidth]{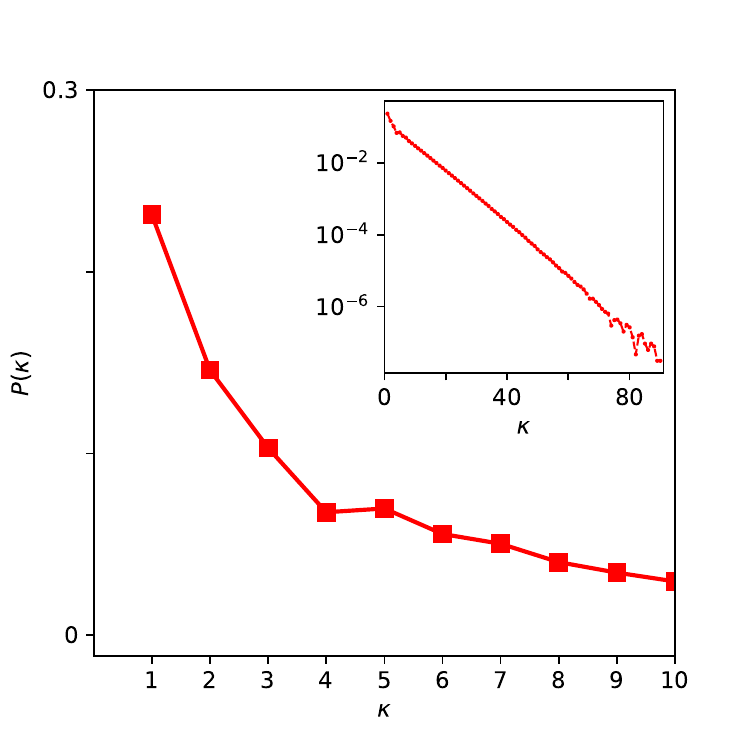}
  \includegraphics[width=0.55\textwidth]{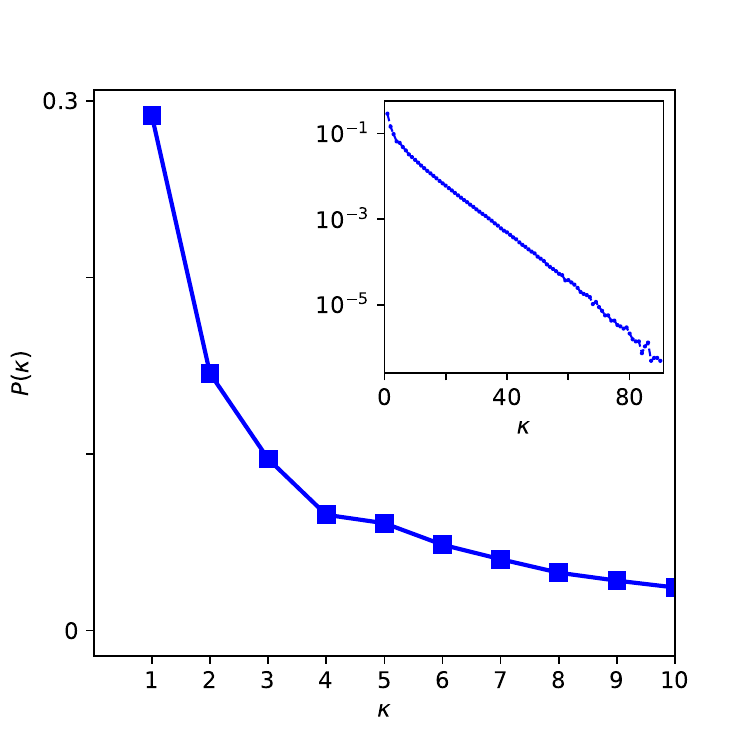}}
  \caption{Degree distribution of random generalised $1$-vertex triangulations of the $3$-dimensional sphere. Top left: fraction of triangulations that have $N(\kappa=1)$ edges with degree $\kappa=1$ for different sizes $n$. Top right: average and standard deviation of the histograms in the left panel (after dividing $N(\kappa)$ by $n-1$) as a function of $n$. Bottom left: the degree distribution $P(\kappa)$ -- estimated at $n=100$ -- showing a non-monotonic behaviour for small $\kappa$ and an exponential decay for large $\kappa$ (straight line in the inset). Bottom right: for comparison, the degree distribution of vertex degrees in generalised triangulations of the $2$-sphere.
  The $\sim 1/\sqrt{n}$ decay of the standard deviation (top right) suggests that almost all large $1$-vertex triangulations of the $3$-sphere have the same $P(\kappa)$. Results measured over the triangulations sampled in \Cref{fig:R,fig:R3}.} \label{fig:degree}
\end{figure}

\subsection{Properties of triangulations}
\label{ssec:vertexdegrees}

A major advantage of our method is that it is not restricted to the counting of triangulations (the question addressed in \Cref{ssec:numberoftriangulations}) or to the finding of specific triangulations (optimisation). Instead, it obtains an unbiased sample of triangulations that can be used to study the prevalence of any triangulation property of interest. In this section, we illustrate this point by analysing a simple property of a generalised triangulation of the $3$-sphere: the sequence of edge degrees. 

The degree of an edge $e$ in a $3$-dimensional triangulation $\tri$, denoted by $\kappa = deg_\tri(e)$, is defined as the number of tetrahedra it is a part of, counted with multiplicity. We start by counting how many degree one ($\kappa=1$) edges we have in $1$-vertex $n$-tetrahedra triangulations of the $3$-sphere $\tri$. As triangulations $\tri$ with a fixed $n$ are sampled with equal frequency, it is enough to count how many triangulations in a given sample have $N(\kappa=1)$ edges with degree $\kappa=1$. This is shown in \Cref{fig:degree} on the left for different values of $n$. Naturally, the number of degree one edges grows with the size of the triangulation $n$ (and hence its number of edges $f_1=n+1$, see \Cref{eq:dim3}). Significantly, the distribution over different triangulations $\tri$ is peaked with a well-defined width. 
We thus characterise the distribution by the mean and standard deviation $\sigma$, after normalising $N(\kappa)$ by the total number of edges $f_1=n+1$. The results for these quantities for different $\kappa$ and $n$ are shown in \Cref{fig:degree} (middle) and indicate a convergence of the means to non-zero values and of the standard deviations to zero (for large $n$). In other words, choosing random triangulations of increasing size $n\rightarrow \infty$, we have that $P(\kappa=1) \approx 23\%$ of the edges have degree one, $P(\kappa=2)=14.5\%$ of the edges have degree two, etc.

In \Cref{fig:degree} on the right we show the degree distribution $P(\kappa)$, which contains several interesting features. First, it is not peaked around the average degree $\langle \kappa \rangle = \frac{6n}{n+1}$ ($6n$ edges for each tetrahedra) but instead decays from $\kappa=1$ (i.e., the most common degree is $\kappa=1$). Second, it is not a monotonic decay function of $\kappa$ but instead there exists at least one local maximum for degree five edges. Third, for large $\kappa$, $P(\kappa)$ decays exponentially with a rate $0.165 \pm 0.06 \approx 1/6$. This rate of decay is directly affected by the erratic behaviour at small $\kappa$ because $\langle \kappa \rangle \equiv \sum_{\kappa=1}^{\infty} \kappa P(\kappa) = 6n/(n+1) \to 6$ (for large $n$).

\section{Conclusions and future directions}

We have introduced and applied a MCMC method to experimentally study triangulations of manifolds. 
The numerical experiments described in this paper can be grouped into two groups. The first was conducted on triangulations of surfaces, largely confirming existing results in the literature. The second investigated triangulations of $3$-manifolds, finding that the number of isomorphism types of triangulations with one vertex and $n$ tetrahedra of the 3-dimensional sphere is consistent with a singly exponential growth in $n$ (if the number of isomorphism classes is super-exponential, then this behaviour is not relevant in the class of $1$-vertex triangulations up to $n=250$).

Our experimental results lead to the following conjectures: (i) the exponent of exponential growth of the number of generalised triangulations $|\Omega(n)|$ of closed orientable surfaces of any genus $g$ is $2.827\pm 0.003$; (ii) the exponent of exponential growth of $|\Omega(n)|$ of generalised $1$-vertex triangulations of the $3$-sphere is $1.5085 \pm 0.0012$; and (iii) the limit of the distribution of the edge degree sequence of any sequence of $1$-vertex generalised $3$-sphere triangulations of diverging size follows an exponential decay for large degrees, but shows non-trivial fluctuations for small degrees.

Here we showed results obtained using a particular choice of MCMC (namely, with acceptance $A=1$). Our general framework can be adapted or extended to different target distributions (e.g., depending on additional properties of $\tri$ beyond $n$ to sample triangulations of particular interest) and different MCMC methods (e.g., beyond Metropolis Hastings~\cite{landau_guide_2014,Fischer2015,prd}). A key improvement on the numerical performance for large $n$ would be achieved avoiding the enumeration of isomorphism types at each step, which would allow the sampling of triangulations at significantly larger $n$ as done in~\cite{prd}.
Another relevant extension would be to apply similar methods beyond the class of triangulations of manifolds (e.g., to more general hypergraphs~\cite{Young2017}).

\bibliographystyle{plain}
\bibliography{references}

\appendix

\section{Appendices}

\subsection{Computation of the acceptance in \Cref{ssec:acceptance}}\label{app.acceptance}

We start our computation of the acceptance in \Cref{eq.A} by computing the ratio of proposals. From \Cref{eq.gi}, with the choice in \Cref{eq.alphatilde}, we obtain for going up ($\delta = n'-n >0$)
 
\begin{equation}\label{eq.gratioU}
  \dfrac{g(\tri'\mapsto\tri)}{g(\tri \mapsto \tri')} = \frac{\tilde{\alpha}(n+\delta)}{\alpha(n)}=\frac{1}{r}\frac{1-\alpha(n+\delta)}{\alpha(n)},
\end{equation}
and for for going down ($\delta = n'-n <0$)
\begin{equation}\label{eq.gratioD}
  \dfrac{g(\tri'\mapsto\tri)}{g(\tri \mapsto \tri')} = \frac{\alpha(n+\delta)}{\tilde{\alpha}(n)} =  r \frac{\alpha(n+\delta)}{1-\alpha(n)}.
\end{equation}

Introducing \Cref{eq.gratioU,eq.gamma,eq.PratioS} in \Cref{eq.A}, we obtain the acceptance for going up ($\delta = n'-n>0$) as
\begin{equation}\label{eq.Aspecificup}
A_{up} = e^{-\beta(2n\delta+\delta^2)} \dfrac{1}{r} \dfrac{1-e^{-\gamma(n+\delta)}}{e^{-\gamma n}} =\frac{e^{-\beta \delta^2}}{r} (1-e^{-\gamma(n+\delta)}) \; e^{+n(\gamma-2\beta\delta)}.
\end{equation}
Similarly,  using \Cref{eq.gratioD,eq.gamma,eq.PratioS} in \Cref{eq.A}, we obtain that the acceptance for going down ($\delta = n'-n<0$) as
\begin{equation}\label{eq.Aspecificdown}
A_{down} = e^{-\beta(2n\delta+\delta^2)} r  \dfrac{e^{-\gamma(n+\delta)}}{1-e^{-\gamma n}} =\frac{r e^{-\beta \delta^2}}{(1-e^{-\gamma n})} \; e^{-n(\gamma+2\beta\delta)}.
\end{equation}
We retrieve \Cref{eq.Aspecificup} noting that $\delta \in [-d,d]$ and considering $n\gg \gamma$ so that $1-e^{-\gamma n} \approx 1$.

\subsection{Numerical Implementation}\label{app.numerics}

{\bf Details of the implementation:}
We use the support for bi-stellar moves built into the low-dimensional topology software {\it Regina} \cite{regina} to implement our method. The algorithm is coded in {\it python3}. The code for the $2$- and $3$-dimensional method can be found in \cite{gitHub}. Our experiments complete the $10\ 000\ 000$ steps in between a couple of minutes and more than a month, depending on the dimension and the value of $\gamma$. The bottleneck for our method is to compute the isomorphism types of the neighbours of our current state triangulation (see \Cref{ssec:triangulations} and the description of the isomorphism signature therein for details on how this is achieved). There is a significant potential to speed up the method through a more sophisticated implementation of enumerating isomorphism types among the neighbour triangulations; or to ignore this step altogether for triangulations of sufficient size, as done in \cite{prd} for a different MCMC method.

\medskip
\noindent
{\bf Data collection:} Given a seed triangulation $\tri$ of $\manifold$, where $\manifold$ is either a surface or a $3$-manifold, we run \Cref{algo:dim2,algo:dim3} for $10\ 000\ 000$ steps with parameter $\gamma = 1/k$, $k \in \{ 2, 3, \ldots , \ell \}$. A larger value of $k$ (a smaller value of $\gamma$) translates to a smaller penalty for choosing a bi-stellar move to a larger triangulation. Hence, for large $\gamma$, small triangulations are sampled and vice versa. We add every $100$th triangulation we see on our walk through the Pachner graph to our sample. 

The result is a list of $100\ 000$ triangulations per run of the experiment. Typically, we repeat every such experiment $20$ times using a different seed for the random generations. We consider the results obtained in each run to be independent estimation of the quantities of interest and thus use the average and standard deviation of the mean (over the $20$ runs) to obtain the reported estimator and corresponding confidence interval.

The parameters of the experiment are chosen such that the probability of (an isomorphism type of) a triangulation to be included in the sample is only distorted by the penalty to go to a larger triangulation, and otherwise approximately uniform. The sample can then be analysed directly, or a subsample can be taken to test the convergence and the general behaviour of the method.

\subsection{Reweighting}\label{app.reweighting}

The Metropolis-Hastings MCMC ensures that for $t\rightarrow \infty$ the probability of sampling a triangulation $\tri$ converges to $P(\tri)$, which in our case depends only on the number of triangles $n(\tri)$. Often one is interested in estimating quantities attributing a different probability (weight) to each triangulation, most commonly an equal weight $P(\tri)=$ constant. One example is the estimation of the (relative) number of triangulations of different types (e.g., different $n$'s). In these cases, an estimation based solely on the sampled data would be biased by $P(\tri)$ (or, equivalently, biased by the uneven proposal and acceptance steps of our random walk). An accurate estimation is obtained by undoing this bias and reweighting the samples.

The key idea of reweighting methods~\cite{landau_guide_2014} is to divide each of the counts (based on $\tri$) by their corresponding sampled probability ($P(\tri)$). The final estimation is typically a relative or reweighted probability (normalisation can be imposed at the end) and therefore knowledge of the ratio of $P(\tri)$ as known here is enough. For instance, consider the experiments in  \Cref{ssec:higherGenus}. Each MCMC run allows us to estimate the ratio of the proportions of $(n+2)$-triangle triangulations over $n$-triangle triangulations. Let the {\it sampled} ratio be $q \in \mathbb{R}$, i.e., $q$ is the number of sampled $\tri$ with $n+2$ triangles divided by the number of sampled $\tri$ with $n$ triangles. The rate of the number of such $\tri$
$$R(n) = \frac{|\Omega(n+2)|}{|\Omega(n)|}$$ 
is then estimated for different $n$ as
$$ \hat{R}(n) = q \cdot \frac{P(n)}{P(n+2)} = q \cdot \frac{1-e^{-\gamma (n+2)}}{2 e^{-\gamma n}},$$
where we used Equation~(\ref{eqPratioS2}) with $r=2$ and $\delta=2$.
For each of the $20$ runs, these estimates are then combined as a weighted average across all parameters, requiring a lower threshold on relative frequencies for each data set of $0.01$ for estimates to be included to eliminate outliers.  Multiple values of $\gamma$ are combined (using a weighted average) to estimate each value of $n$ because every set of samples coming from a sequence of moves with fixed $\gamma$ gives samples with a range of sizes. Finally, the average over all $20$ runs is computed.  

\subsection{Enumeration}\label{sec:enumeration}

\Cref{algo:dim2} can be adapted to yield an ad-hoc enumeration procedure for $n$-triangle generalised triangulations of surfaces of a given genus $g$ as follows: 
\begin{itemize}
  \item Start with a seed $n$-triangle triangulation of a surface of genus $g$. This is obtained from a minimal $4g-2$-triangle triangulation of that surface followed by an appropriate number $0$-moves.
  \item Search through the Pachner graph of $n$-triangle generalised triangulations of this surface by deterministically performing all $1$-moves on the seed triangulation and then iterate in a bread-first fashion.
  \item Record all isomorphism types encountered in this search until the list of isomorphism types becomes stationary.
\end{itemize}
By construction, this procedure delivers a complete list of all $n$-triangle genus $g$ surface triangulations. This enumeration procedure is, while certainly not very efficient, very convenient to implement and fast enough to produce gigabytes of triangulations in moderate time frames.

\end{document}